\documentclass[12pt,a4paper]{article}
\usepackage{amssymb}
\usepackage[utf8]{inputenc}
\usepackage[english]{babel}
\usepackage{color}
\usepackage{amsmath}
\usepackage{enumerate}
\usepackage{amssymb,enumerate}

\newtheorem{Thm}{Theorem}[section]
\newtheorem{Def}[Thm]{Definition}
\newtheorem{Prop}[Thm]{Proposition}
\newtheorem{Lem}[Thm]{Lemma}
\newtheorem{Cor}[Thm]{Corollary}
\newtheorem{Rem}[Thm]{Remark}

\title{$L^p$-distortion and $p$-spectral gap of finite graphs}
\author{Pierre-Nicolas JOLISSAINT\footnote{Supported by Swiss SNF project 20-137696.}$\;$ and Alain VALETTE}

\begin{document}

\maketitle

\begin{abstract} We give a lower bound for the $L^p$-distortion $c_p(X)$ of finite graphs $X$, depending on the first eigenvalue $\lambda_1^{(p)}(X)$ of the $p$-Laplacian and the maximal displacement of permutations of vertices. For a $k$-regular vertex-transitive graph it takes the form $$c_p(X)^{p}\geq \frac{diam(X)^{p}\lambda_{1}^{(p)}(X)}{2^{p-1}k}.$$ This bound is optimal for expander families and, for $p=2$, it gives the exact value for cycles and hypercubes. As new applications we give non-trivial lower bounds for the $L^2$-distortion for families of Cayley graphs of the finite lamplighter groups $C_{2}\wr C_{n}^{d}$ ($d\geq 2$ fixed), and for a family of Cayley graphs of $SL_n(q)$ ($q$ fixed, $n\geq 2$) with respect to a standard two-element generating set. An application to the $L^2$-compression of certain box spaces is also given.
\end{abstract}

\section{Introduction}
Let $(X,d)$ and $(Y,\delta)$ be two metric spaces. Let $F : X \rightarrow Y$ be an imbedding of $X$ into $Y$. We define the \textit{distortion} of $F$ as
$$
dist(F)=\sup_{x,y\in X, x\neq y} \frac{\delta(F(x),F(y))}{d(x,y)} \cdot \sup_{x,y\in X, x\neq y} \frac{d(x,y)}{\delta(F(x),F(y))},
$$
where the first supremum is the Lipschitz constant $\|F\|_{Lip}$ of $F$, and the second supremum is the Lipschitz constant $\|F^{-1}\|_{Lip}$ of $F^{-1}$. As we will only consider the case where $X$ is finite, supremum can be changed into maximum. The least distortion with which $X$ can be embedded into $Y$ is denoted $c_{Y}(X)$, namely
$$
c_{Y}(X):=\inf\{ dist(F) : \ F : X \hookrightarrow Y \}.
$$
As target space, we will consider only $L^{p}=L^p([0,1])$. In this case, we write $c_{p}(X)=c_{L^{p}}(X)$. The quantity $c_{2}(X)$ is also known as the Euclidean distortion of $X$. As source space, we will take the underlying metric space of a finite, connected graph $X=(V,E)$, where $d$ is then the graph metric. Note that, denoting by $diam(X)$ the diameter of $X$, we have $c_p(X)\leq diam(X)$, as shown by the embedding $F:V\rightarrow\ell^p(V):x\mapsto\delta_x$. It is a fundamental result of Bourgain \cite{Bou} that\footnote{In this paper, Landau's notations $O,\,\Omega,\,\Theta$ will be used freely.} : $c_p(X)=O(\log |V|)$.

Our aim in this paper is to obtain lower bounds for the distortion $c_{p}$ of finite graphs. To state our results, we introduce two invariants of graphs. For $1<p< \infty$, the $p$-Laplacian $\Delta_{p} : \ell^{p}(V) \rightarrow \ell^{p}(V)$ is an operator defined by the formula
$$
\Delta_{p}f(x)=\sum_{x\sim y} \left( f(x)-f(y)\right)^{\left[ p\right]},
$$
($f\in\ell^{p}(V), x\in V$), where $a^{\left[p\right]}=|a|^{p-1}sign(a)$ and $\sim$ denotes the adjacency relation on $V$. It is worth noting that for $p=2$, it corresponds to the standard linear discrete Laplacian. We say that $\lambda$ is an eigenvalue of $\Delta_{p}$ if we can find $f\in\ell^{p}(V)$ such that $\Delta_{p}f=\lambda f^{\left[p\right]}$. For $1\leq p < \infty$, we define the \textit{$p$-spectral gap} of $X$ by
$$
\lambda_{1}^{(p)}(X):=\inf\left\{\frac{\frac{1}{2}\sum_{x\in V}\sum_{y:\;y\sim x}|f(x)-f(y)|^{p}}{\inf_{\alpha\in\mathbb{R}} \sum_{x\in V}|f(x)-\alpha|^{p}}\right\},
$$
where the infimum is taken over all $f\in\ell^{p}(V)$ such that $f$ is not constant. When, $p\neq 1$, it is known that the $p$-spectral gap is the smallest positive eigenvalue of $\Delta_{p}$ (see \cite{GN}).

For $\alpha$ a permutation of the vertex set $V$ (not necessarily a graph automorphism!), we introduce the {\it displacement} of $\alpha$:
$$\rho(\alpha)=\min_{x\in V} d(\alpha(v),v);$$
then the {\it maximal displacement} of $X$ is $D(X):=\max_{\alpha\in Sym(V)}\rho(\alpha)$. (Note that this definition makes sense for every finite metric space).

Our main result is:

\begin{Thm}\label{generalization}
Let $X$ be a finite, connected graph of average degree $k$. Then
$$
D(X) \left( \frac{\lambda_{1}^{(p)}(X)}{k \ 2^{p-1}}\right)^{\frac{1}{p}} \leq c_{p}(X),
$$
for $1\leq p<\infty$.
\end{Thm}

For vertex-transitive graphs, this takes the form:

\begin{Cor}\label{generalizationCayley}
Let $X$ be a finite, connected, vertex-transitive graph. Then for $1\leq p<\infty$:
$$
diam(X) \left( \frac{\lambda_{1}^{(p)}(X)}{k \ 2^{p-1}}\right)^{\frac{1}{p}}\leq c_{p}(X),
$$
where $k$ is the degree of each vertex.
\end{Cor}

Theorem \ref{generalization} and Corollary \ref{generalizationCayley} allow us to give unified proofs of some known results:

\begin{enumerate}

\item[1)] (see \cite{LLR},\cite{matou}\footnote{Expanders were used by Linial-London-Rabinovich \cite{LLR} for $p=2$, and by Matou\v sek \cite{matou} for arbitrary $p\geq 1$, to show that Bourgain's upper bound on $c_p$ is optimal for every $p$.}) Let $(X_{n})_{n\geq 1}$ be a family of expanders with bounded degree. For every $p \geq 1$, we have  $c_p(X_{n})=\Omega(\log |X_{n}|).$
\item[2)] (Linial-Magen \cite{LM}) For even $n$: the cycle $C_n$ satisfies $c_2(C_n)=\frac{n}{2}\sin\frac{\pi}{n}$.
\item[3)] (Enflo \cite{Enf}) The $d$-dimensional hypercube $H_d$ satisfies $c_2(H_d)=\sqrt{d}$.

\end{enumerate}

\medskip

As a first application, we consider lamplighter groups over discrete tori. Recall that, if $G$ is a finite group, the lamplighter group of $G$ is the wreath product $C_2\wr G$, i.e.\ the semi-direct product of the additive group of all subsets of $G$ (endowed with symmetric difference) with $G$ acting by shifting indices. Take $G=C_n^d$ and denote by $\{\pm e_{j} \ : \ 1\leq j \leq d\}$ the standard symmetric generating set for $C_{n}^{d}$, and denote by $W_{n}^{d}$ the Cayley graph of the lamplighter group $C_2\wr C_{n}^{d}$, with respect to the generating set 
$$
S=\{(\{0\},0)\} \cup\{(\emptyset,\pm e_j):1\leq j\leq d\}.
$$
(so that $W_n^d$ is $(2d+1)$-regular). We will prove the following:

\begin{Prop}\label{wreath} $c_2(W_{n}^{d})=\left\{
\begin{array}{ll}
\Omega(\frac{n}{\sqrt{\log(n)}}), & \textrm{for}\ d=2,\\
\Omega(n^{\frac{d}{2}}), & \textrm{for}\ d\geq 3.\\
\end{array}\right.$
\end{Prop}

However, the method we will use does not give a good estimate for the case $d=1$ as we will see in section $5$.

\medskip
As a second application, let $q$ be a fixed prime, and let $Y_n$ be the Cayley graph of $SL_n(q)$ (where $n\geq 2$) with respect to the following set of 4 generators: $S_n=\{A_n^{\pm 1},B_n^{\pm 1}\}$ and 
$$A_n=\left(\begin{array}{ccccc}
1 & 1 & & & \\
& 1 & & & \\
& & 1 & & \\
& & & \ddots & \\
& & & & 1
\end{array}\right);\;\;B_n=\left(\begin{array}{cclll}
0 & 1 & & & \\
& 0 & 1 & & \\
& & 0 & \ddots & \\
& & & \ddots & 1 \\
(-1)^{n-1} & & & & 0
\end{array}\right).
$$

\begin{Prop}\label{slnp} $c_2(Y_n)=\Omega(n^{1/2})=\Omega((\log |Y_n|)^{1/4})$.
\end{Prop}

The interest of the family $(Y_n)_{n\geq 2}$ comes from the fact that it is known NOT to be an expander family: see Proposition 3.3.3 in \cite{Lub}. 

\medskip
The paper is organized as follows: Theorem \ref{generalization} is proved in section 2, and Corollary \ref{generalizationCayley} in section 3, where estimates on the maximal displacement are also given. Various applications are given in section 4, which also presents examples where the inequality in Corollary \ref{generalizationCayley} is {\it not} sharp. The section 5 contains a discussion of  other published results similar to our Theorem \ref{generalization}, and a comparison of the corresponding inequalities. Finally, in section 6, we use the estimates from section 4 to give non-trivial upper bounds for the compression exponent of some infinite metric spaces obtained as disjoint union of finite graphs.

\medskip

{\bf Acknowledgements:} We thank R. Bacher, B. Colbois, A. Gournay, A. Lubotzky, R. Lyons and Y. Peres for useful exchanges, and comments on the first draft. Some of the results of this paper were included in the pre-book \cite{LyPe}.

\pagebreak

\section{Proof of Theorem \ref{generalization}}

We start with an easy lemma. 
\begin{Lem}\label{lem2}
Let $X=(V,E)$ be a finite, connected graph and let $1<p<\infty$.
\begin{enumerate}
	
	\item Let $\alpha$ be any permutation of $V$. For $F : V \rightarrow \ell^{p}(\mathbb{N})$ :
	$$
	\sum_{x\in V} \|F(x)-F(\alpha (x))\|_{p}^{p} \leq 2^p \sum_{x\in V} \|F(x)\|_{p}^{p}.
	$$
	
	\item Fix an arbitrary orientation on the edges. Then, for every $F : V \rightarrow \ell^{p}(\mathbb{N})$, there exists $G : V \rightarrow \ell^{p}(\mathbb{N})$ such that $dist(G)=dist(F)$ and
	$$
	\sum_{x\in V} \|G(x)\|_{p}^{p} \leq \frac{1}{\lambda_{1}^{(p)}(X)} \sum_{e\in E} \|G(e^{+})-G(e^{-})\|_{p}^{p}.
	$$
\end{enumerate}
\end{Lem}

{\bf Proof:} 1) Define a linear operator $T$ on $\ell^{p}(V, \ell^{p}(\mathbb{N}))$ by setting $(TF)(x):=F(\alpha (x))$.  Clearly, $\|T\|=1$. Then, in the formula to be proved, the LHS is $\|(I-T)F\|_{p}^{p}$. Hence, the result immediately follows from the fact that the operator norm of $T-I$ is at most $2$, by the triangle inequality. \\
2) We proceed as in the proof of Theorem $3$ in \cite{GN}. Let $\{u_{n}\}_{n\in\mathbb{N}}$ be the standard basis vectors in $\ell^{p}(\mathbb{N})$. 
Write $F(x)=\sum_{n\in\mathbb{N}}F_{n}(x)u_{n}$, for all $x\in V$; using the fact that $\ell^{p}(V)$ is uniformly convex for $p>1$, we denote by $\alpha_{n}\in\mathbb{R}$ the projection of $F_n$ on the subspace of constant functions in $\ell^p(V)$. It satisfies:
$$
\inf_{\alpha\in\mathbb{R}} \|F_{n}-\alpha\|_p = \|F_{n}-\alpha_{n}\|_p.
$$
By the proof of Theorem  3 in \cite{GN}, the sum $w:=\sum_{n\in\mathbb{N}} \alpha_{n}u_{n}$ belongs to $\ell^{p}(\mathbb{N})$. 

Defining $G(x):=F(x)-w$, so that $G_{n}(x)=F_{n}(x)-\alpha_{n}$, we have $dist(G)=dist(F)$.
Recalling the definition of $\lambda_{1}^{(p)}(X)$, we have for every $n$:
$$
\sum_{x\in V} |G_{n}(x)|^{p} \leq \frac{1}{\lambda_{1}^{(p)}(X)}\sum_{e\in E} |G_{n}(e^{+})-G_{n}(e^{-})|^{p}.
$$
Taking the sum over $n$, we get the result.
\hfill $\square$ \\

\medskip

Let $k$ be the average degree of $X$. Combining both statements of lemma \ref{lem2} with the fact that $|E|=\frac{k|V|}{2}$, we deduce the following Poincar\'e-type inequality:
\begin{Prop}\label{Poincare}
Let $1<p<\infty$ and let $X=(V,E)$ be a finite, connected graph with average degree $k$. For any permutation $\alpha$ of $V$ and any embedding $G : V \rightarrow \ell^{p}(\mathbb{N})$ as in lemma \ref{lem2},  we have:
$$
\frac{1}{|V|2^{p}} \sum_{x\in V} \|G(x)-G(\alpha (x))\|_{p}^{p} \leq \frac{k}{2|E|\lambda_{1}^{(p)}(X)} \sum_{e\in E} \|G(e^{+})-G(e^{-})\|_{p}^{p}.
$$
\end{Prop}
\hfill $\square$



\begin{Prop}\label{Prop.avgVvsE}
Let $X=(V,E)$ be a finite connected graph with average degree $k$. For any permutation $\alpha$ of $V$ and any embedding $G : V \rightarrow \ell^{p}(\mathbb{N})$ as in lemma \ref{lem2}, we have:
$$
\rho (\alpha) \left( \frac{\lambda_{1}^{\left(p\right)}(X)}{k \ 2^{p-1}}\right)^{\frac{1}{p}} \leq dist(G).
$$

\end{Prop}

{\bf Proof:} Clearly, we may assume that $\alpha$ has no fixed point. Then:
$$
\frac{1}{\|G^{-1}\|_{Lip}^{p}}=\min_{x\neq y} \frac{\|G(x)-G(y)\|_{p}^{p}}{d(x,y)^{p}} \leq \min_{x\in V}  \frac{\|G(x)-G(\alpha (x))\|_{p}^{p}}{d(x,\alpha (x))^{p}}
$$
$$
\leq \frac{1}{\rho(\alpha)^{p}} \min_{x\in V} \|G(x)-G(\alpha (x))\|_{p}^{p} \leq \frac{1}{\rho(\alpha)^{p}|V|} \sum_{x\in V} \|G(x)-G(\alpha (x))\|_{p}^{p}
$$
$$
\leq \frac{2^{p-1}k}{\lambda_{1}^{(p)}(X)\rho(\alpha)^{p}|E|}\sum_{e\in E} \|G(e^+)-G(e^-)\|_{p}^{p} \; \; \mbox{(by Proposition \ref{Poincare})}
$$
$$
\leq \frac{2^{p-1}k}{\lambda_{1}^{(p)}(X)\rho(\alpha)^{p}} \ \max_{x\sim y} \|G(x)-G(y)\|_{p}^{p}= \frac{2^{p-1}k}{\lambda_{1}^{(p)}(X)\rho(\alpha)^{p}} \|G\|_{Lip}^{p},
$$
where the last equality comes from the fact that the above maximum is attained for adjacent points in the graph (see for instance Claim $3.2$ in \cite{LM}). Re-arranging and taking $p$-th roots, we get the result.
 \hfill $\square$
 
 \medskip
 
 {\bf Proof of Theorem \ref{generalization}}: Since $\ell^p$ embeds isometrically in $L^p$, we clearly have $c_p(X)\leq c_{\ell^p}(X)$. Actually $c_p(X)=c_{\ell^p}(X)$, since for every map $F:V\rightarrow L^p$ and every $\varepsilon>0$, we can find a finite measurable partition $[0,1]=\bigcup_{j=1}^k \Omega_j$ and, for each $x\in V$, a step function $H(x)$ which is constant on each $\Omega_j$, such that $\|F(x)-H(x)\|_p <\varepsilon$ for $x\in V$. Denoting by $m$ the Lebesgue measure on $[0,1]$, the embedding $G:V\rightarrow \ell^p\{1,...,k\}:x\mapsto (H(x)|_{\Omega_j}m(\Omega_j)^{1/p})_{1\leq j\leq k}$ then satisfies $\|G(x)-G(y)\|=\|H(x)-H(y)\|_p$ for every $x,y\in V$, hence the distortion of $G$ is $\delta(\varepsilon)$-close to the one of $F$, where $\delta(\epsilon)\rightarrow 0$ for $\varepsilon\rightarrow 0$. Now, Theorem \ref{generalization} for embeddings $V\rightarrow \ell^p$, where $1<p<\infty$, immediately follows from Proposition \ref{Prop.avgVvsE}. Finally, a straightforward continuity argument allows us to cover the case $p=1$.
\hfill$\square$

\section{Estimates on the maximal displacement}

From the definition of the invariant $D(X)$, we have $D(X)\leq diam(X)$. The equality holds if and only if the graph $X$ admits an {\it antipodal map}, i.e. a permutation $\alpha$ of the vertices such that $d(x,\alpha(x))=diam(X)$ for every $x\in V$. 

The existence of an antipodal map is a fairly strong condition. Recall that the {\it radius} of $X$ is $\min_{x\in V}\max_{y\in V}d(x,y)$, so that the existence of an antipodal map implies that the radius is equal to the diameter of $X$. The converse is false however, a counter-example was provided by G. Paseman. A necessary and sufficient condition for $X$ to admit an antipodal map was provided by R. Bacher: for $S\subset V$, set  ${\cal A}(S)=\{v\in V: \exists w\in S, d(v,w)=diam(X)\}$; the graph $X$ admits an antipodal map if and only if $|{\cal A}(S)|\geq|S|$ for every $S\subset V$. For all this, see \cite{MO}.

The proof of Corollary \ref{generalizationCayley} follows immediately from Theorem \ref{generalization} and the next lemma:

\begin{Lem}\label{Cayley} Finite, connected, vertex-transitive graphs admit antipodal maps.
\end{Lem}

{\bf Proof:} For $S$ a finite subset of the vertex set of some graph $Y$, denote by $\Gamma(S)$ the set of vertices adjacent to at least one vertex of $S$. It is classical that, if $Y$ is a regular graph, then the inequality $|\Gamma(S)|\geq |S|$ holds\footnote{Recall the easy argument: assuming that $Y$ is $k$-regular, count in two ways the edges joining $S$ to $\Gamma(S)$; as edges emanating from $S$, there are $k|S|$ of them; as edges entering $\Gamma(S)$, there are at most $k|\Gamma(S)|$ of them.}.

Now, let $X=(V,E)$ be a finite, connected, vertex-transitive graph. Define the {\it antipodal graph} $X^a$ as the graph with vertex set $V$, with $x$ adjacent to $y$ whenever the distance between $x$ and $y$ in $X$, is equal to $diam(X)$. By vertex-transitivity of $X$, the graph $X^a$ is regular. Now observe that, for $S\subset V$, the set $\Gamma(S)$ in $X^a$ is exactly the set ${\cal A}(S)$ defined above. By regularity of $X^a$ and the observation beginning the proof, we therefore have $|{\cal A}(S)|\geq|S|$ for every $S\subset V$, and Bacher's result applies.
 \hfill$\square$ \\

\begin{Rem} For Cayley graphs, there is a direct proof of the existence of antipodal maps. Indeed, let $G$ be a finite group, and let $X$ be a Cayley graph of $G$ with respect to some symmetric, generating set $S$; use right multiplications by generators to define $X$, so that the distance $d$ is left-invariant. Let $g\in G$ be any element of maximal word length with respect to $S$. Then $\alpha(x)=xg$ (right multiplication by $g$) is an antipodal map. 
\end{Rem}



For arbitrary graphs, we have:

\begin{Prop}\label{displace} For finite, connected graphs $X$ with maximal degree $k\geq 3$:
$$D(X)=\Omega(\log|X|).$$
\end{Prop}

{\bf Proof:} For a positive integer $r>0$, the number of vertices in $X$ at distance at most $r$ from a given vertex, is at most the number of vertices in the ball of radius $r$ in the $k$-regular tree, i.e.
$$1 + k + k(k-1) + k(k-1)^2 + ... + k(k-1)^{r-1} = \frac{k(k-1)^r-2}{k-2}.$$
For $r=[\log_{k-1}(\frac{|V|}{6})]$, we have $\frac{k(k-1)^r-2}{k-2}<\frac{|V|}{2}$. Let $Y$ be the graph with the same vertex set $V$ as $X$, where two vertices are adjacent if their distance in $X$ is at least $\log_{k-1}(\frac{|V|}{6})$. The preceding computation shows that, in the graph $Y$, every vertex has degree at least $\frac{|V|}{2}$. By G.A. Dirac's theorem (see e.g. Theorem 2 in Chapter IV of \cite{Bol}), $Y$ admits a Hamiltonian circuit. Let $\alpha\in Sym(V)$ be the cyclic permutation of $V$ defined by this Hamiltonian circuit. Then $\rho(\alpha)\geq\log_{k-1}(\frac{|V|}{6})$, which concludes the proof.
\hfill$\square$

\medskip

The following inequality due to Alon and Milman (see Theorem 2.7 in \cite{AM}) shows that Proposition \ref{displace} is essentially the best possible. For any connected graph $X=(V,E)$ with degree bounded by $k$, we have
$$
diam(X)\leq 2 \sqrt{\frac{2k}{\lambda_{1}^{(2)}(X)}} \ \log_{2}|X|,
$$

\medskip


We now observe that, for families of non-vertex-transitive $k$-regular graphs, the maximal displacement can be much smaller than the diameter (compare with lemma \ref{Cayley}). We thank the referee of a previous version of the paper for suggesting this construction.

\begin{Prop}\label{stitch} Let $f:\mathbb{N}\rightarrow\mathbb{N}$ be a function such that $f(n)=\Omega(n)$ and $f(n)= o(8^n)$. There exists a family $(X_n)_{n\geq 1}$ of 3-regular graphs such that: 
\begin{enumerate}
\item[a)] $|X_n|=\Theta(8^n)$;
\item[b)] $diam(X_n)=\Theta(f(n))$;
\item[c)] $D(X_n)=\Theta(n)$.
\end{enumerate}

\end{Prop}

{\bf Proof:} Let $(Y_n)_{n\geq 1}$ be a family of 3-regular graphs with $|Y_n|=\Theta(8^n)$ and $diam(Y_n)=\Theta(n)$ (such a family is constructed e.g. in Theorem 5.13 of Morgenstern \cite{Mor}). Let $Z_n$ be the product of the cycle $C_{2f(n)}$ with the one-edge graph (so that $Z_n$ is 3-regular on $4f(n)$ vertices). Let $\{y_1,y_2\}$ (resp. $\{z_1,z_2\}$) be an edge in $Y_n$ (resp. $Z_n$). We ``stitch'' $Y_n$ and $Z_n$ by replacing the edges $\{y_1,y_2\}$ and $\{z_1,z_2\}$ by edges $\{y_1,z_1\}$ and $\{y_2,z_2\}$, and define $X_n$ as the resulting 3-regular graph. Clearly $|X_n|=\Theta(8^n)$. 

Observe that, since every edge in $Z_n$ belongs to some 4-cycle, the distance in $X_n$ between any two vertices in $Y_n$ will differ by at most 5 from the original distance in $Y_n$; and similarly for vertices in $Z_n$. So:
$$f(n) = diam(Z_n)\leq diam(X_n)\leq diam(Y_n) + diam(Z_n) + 5,$$
hence $diam(X_n)=\Theta(f(n))$.

Finally, let $\alpha$ be any permutation of the vertices of $X_n$. Since the overwhelming majority of vertices belongs to $Y_n$, we find a vertex $x$ such that $x$ and $\alpha(x)$ are both in $Y_n$. Then
$$\rho(\alpha)\leq d_{X_n}(x,\alpha(x))\leq d_{Y_n}(x,\alpha(x))+5\leq diam(Y_n)+ 5,$$
hence $D(X_n)=O(n)$. The equivalence $D(X_n)=\Theta(n)$ then follows from Proposition \ref{displace}.
\hfill$\square$

\section{Applications}

We give a series of consequences of Theorem \ref{generalization} and Corollary \ref{generalizationCayley}.

\subsection{Expanders}

\begin{Cor}
(\cite{LLR},\cite{matou}) Let $(X_{n})_{n\geq 1}$ be a family of expanders with bounded degree. For every $p\geq 1$, we have  $c_p(X_{n})=\Omega(\log |X_{n}|).$
\end{Cor}
{\bf Proof:} If $(X_n)_n$ is a family of expanders, then by the $p$-Laplacian version of the Cheeger inequality (see Theorem 3 in \cite{Amg}), the sequence $(\lambda_1^{(p)}(X_n))_n$ is bounded away from $0$. So the result follows straight from Theorem \ref{generalization} together with Proposition \ref{displace}.
\hfill$\square$

\subsection{Cycles}

\begin{Cor}\label{cycle} (Linial-Magen \cite{LM}, 3.1) For $n$ even: $c_2(C_n)=\frac{n}{2}\sin\frac{\pi}{n}.$
\end{Cor}

{\bf Proof:} We apply Corollary \ref{generalizationCayley} with $k=2$, and $D=\frac{n}{2}$, and $\lambda_1^{(2)}(C_n)=4\sin^2\frac{\pi}{n}$ (see Example 1.5 in \cite{Chu}): so $c_2(C_n)\geq \frac{n}{2}\sin\frac{\pi}{n}$. For the converse inequality, it is an easy computation that the embedding of $C_n$ as a regular $n$-gon in $\mathbb{R}^2$, has distortion $\frac{n}{2}\sin\frac{\pi}{n}.$
\hfill$\square$

\subsection{The hypercube $H_d$}

The hypercube $H_d$ is the set of $d$-tuples of $0$'s and $1$'s, endowed with the Hamming distance. It is the Cayley graph of $\mathbb{F}_2^d$ with respect to the standard basis.

\begin{Cor}\label{hypercube} (Enflo \cite{Enf}) $c_2(H_d)=\sqrt{d}$
\end{Cor}

{\bf Proof:} For $H_d$, we have $k=d$, and $diam(H_d)=d$, and $\lambda_1^{(2)}(H_d)=2$ (see Example 1.6 in \cite{Chu} for the latter): so $c_2(H_d)\geq\sqrt{d}$ by Corollary \ref{generalizationCayley}. For the converse inequality, it is easy to see that the canonical embedding of $H_d$ into $\mathbb{R}^d$, has distortion $\sqrt{d}$.
\hfill$\square$

\subsection{Lamplighters over discrete tori}

Once again we apply Corollary \ref{generalizationCayley} in order to prove Proposition \ref{wreath}. Recall that $W_{n}^{d}$ refers to the Cayley graph of the lamplighter group $C_2\wr C_{n}^{d}$ with respect to the standard generating set. Let us define the matrix $M$ on $C_2\wr C_n^d$ given by
$$
M_{\lbrack(f,a),(g,b)\rbrack}=\left\{
\begin{array}{llll}
\frac{1}{4}& \textrm{if}\ (f,a)=(g,b);\\
\frac{1}{4}& \textrm{if}\ a=b \ \textrm{and}\ f=g+\delta_{a} ;\\
\frac{1}{16d}& \textrm{if}\ a=b\pm e_{j}\ \textrm{and}\ f(z)=g(z), \forall z\notin\{a,b\}; \\
0& \textrm{otherwise} .\\
\end{array}\right.
$$
($a,b\in C_{n}^{d}$ and $f,g: C_{n}^{d} \rightarrow \{0,1\}$). Then $M$ is the transition matrix of the lazy random walk on $C_{2} \wr C_{n}^{d}$ analysed by Peres and Revelle in Theorem 1.1 of  \cite{PR}. Using their estimation of the relaxation time of $M$, we deduce that the spectral gap of $M$ behaves as  $\Theta(\frac{1}{n^{d}})$ for $d\geq 3$ and as $\Theta(\frac{1}{n^{2}\log(n)})$ for the case $d=2$. By standard comparison theorems (see e.g. Theorems 3.1 and 3.2 in \cite{Woe}), the Dirichlet forms for $M$ and for the Laplace operator on $W_{n}^{d}$ are bi-Lipschitz equivalent ; moreover the Lipschitz constants do not depend on $n$ (since the comparison can be made on the group $C_2\wr\mathbb{Z}^d$, of which our lamplighters are quotients). So, we find $\lambda_{1}^{(2)}(W_{n}^{2})=\Theta(n^{-2}\log(n)^{-1})$ and $\lambda_{1}^{(2)}(W_{n}^{d})=\Theta(n^{-d})$ for $d\geq 3$. Furthermore, since the diameter of a regular graph is at least logarithmic in the number of vertices, we have $diam(W_{n}^{d})=\Omega(n^{d})$, so we apply Corollary \ref{generalizationCayley} to get:
$$
c_2(W_{n}^{d})=\left\{
\begin{array}{ll}
\Omega(\frac{n}{\sqrt{\log(n)}})& \textrm{for}\ d=2,\\
\Omega(n^{\frac{d}{2}})& \textrm{for}\ d\geq 3.\\
\end{array}\right.
$$  

\hfill$\square$

\subsection{Cayley graphs of $SL_n(q)$}

We now prove Proposition \ref{slnp}. Since $|SL_n(q)|\approx q^{n^2-1}$, we have $diam(Y_n)=\Omega(n^2)$ (actually it is a result by Kassabov and Riley \cite{KaRi} that $diam(Y_n)=\Theta(n^2)$). On the other hand, from Kassabov's estimates for the Kazhdan constant $\kappa(SL_n(\mathbb{Z}),S_n)$ (see \cite{Kas}, and also the Introduction of \cite{KaRi}), we have: $\kappa(SL_n(\mathbb{Z}),S_n)=\Omega(n^{-3/2})$.

If $X$ is a Cayley graph of a finite quotient of a Kazhdan group $G$, with respect to a finite generating set $S\subset G$, then $\lambda_1^{(2)}(X)\geq\frac{\kappa(G,S)^2}{2}$ (see \cite{Lub}, Proposition 3.3.1 and its proof). From this we get: $\sqrt{\lambda_1^{(2)}(Y_n)}=\Omega(n^{-3/2})$ and therefore $c_2(Y_n)=\Omega(n^{1/2})$ by Corollary \ref{generalizationCayley}.
\hfill$\square$

\subsection{The limits of the method}

We give examples of Cayley graphs for which the lower bound of the Euclidean distortion given by Corollary \ref{generalizationCayley} is not tight. 

\subsubsection{Products of cycles}

Let us consider the product of $2$ cycles $C_{n}\times C_{N}$, where $n, N$ are even integers such that $n<N$. It is clear that it corresponds to the Cayley graph of the additive group $\mathbb{Z}/ n\mathbb{Z}\times\mathbb{Z}/ N\mathbb{Z}$ with generating set $S=\{(\pm 1,0),(0,\pm 1)\}$. It is well-known from representation theory of finite abelian groups $G$ that, if $X=\mathcal{G}(G,S)$ is a Cayley graph of $G$ and $S$ is symmetric, then the spectrum of the Laplace operator on $X$ is given by $\{\sum_{s\in S} (1-\chi(s)) \ : \ \chi\in\hat{G}\}$. Since for the product of finite abelian groups $G,H$, we can identify the dual of $G \times H$ as $\{\chi \cdot \eta \ : \ \chi\in\hat{G}, \eta\in\hat{H}  \}$, it is easy to see that $\lambda_{1}(C_{n}\times C_{N})=4\sin^2 \frac{\pi}{N}$. As the diameter is equal to $\frac{n+N}{2}$, we get the lower bound
$$
c_{2}(C_{n}\times C_{N}) \geq \frac{(n+N)\sin\frac{\pi}{N}}{2\sqrt{2}}.
$$
On the other hand, it is known from \cite{LM} that the normalized trivial embedding of $C_{n}\times C_{N}$ into $\mathbb{C}^{2}$ gives the optimal embedding. Namely, defining
$$
\phi : C_{n}\times C_{N} \rightarrow \mathbb{C}^{2} : (k,l) \mapsto \left(\frac{\exp{\frac{2\pi ik}{n}}}{2\sin\frac{\pi}{n}},\frac{\exp{\frac{2\pi il}{N}}}{2\sin\frac{\pi}{N}}\right)
$$   
we have
$$
c_{2}(C_{n}\times C_{N})=dist(\phi).
$$
Since $\|\phi(x)-\phi(y)\| \leq 1$ for every $x,y\in C_{n}\times C_{N}$, we have to estimate
$$
\|\phi^{-1}\|_{Lip}=\max_{k\leq \frac{n}{2}, l\leq \frac{N}{2}} \frac{k+l}{\sqrt{\frac{sin^{2}\frac{\pi k}{n}}{sin^{2}\frac{\pi}{n}}+\frac{sin^{2}\frac{\pi l}{N}}{sin^{2}\frac{\pi}{N}}}}.
$$
By taking $k=\frac{n}{2}$ and $l=\frac{N}{2}$, we get
$$
dist(\phi)\geq \frac{n+N}{2\sqrt{\sin^{-2}\frac{\pi}{n}+\sin^{-2}\frac{\pi}{N}}}. 
$$
Since it is always the case that 
$$
\sqrt{\frac{1}{\sin^{-2}\frac{\pi}{n}+\sin^{-2}\frac{\pi}{N}}} > \frac{\sin \frac{\pi}{N}}{\sqrt{2}},
$$
we conclude that the lower bound given by Corollary \ref{generalizationCayley} is not sharp in this case.

\subsubsection{Lamplighter groups over the discrete circle}

Here we consider the graphs $W_n^1$ associated with the lamplighter groups $C_2\wr C_n$, associated with the generating $S$ described in the Introduction. It is known from \cite{ANV} that $c_2(W_n^1)=\Theta(\sqrt{\log(n)})$. By way of contrast, let us check that $diam(W_n^1)\sqrt{\lambda_1^{(2)}(W_n^1)}=O(1)$. Let us first estimate the spectral gap. For every homomorphism $\chi:C_2\wr C_n\rightarrow\mathbb{C}^\times$, the quantity $\sum_{s\in S}(1-\chi(s))$ is an eigenvalue of the Laplace operator (see the previous example). Let us consider the homomorphism $\chi$ given by $\chi(A,k)=e^{2\pi ik/n}$ (it factors through the epimorphism $C_2\wr C_n\rightarrow C_n$). Here we get $\lambda_1^{(2)}(W^1_n)\leq \sum_{s\in S}(1-\chi(s))= 2-2\cos(2\pi/n)=4\sin^2(\pi/n)$, hence $\lambda_1^{(2)}(W^1_n)=O(\frac{1}{n^2})$. On the other hand, by Theorem 1.2 in \cite{Par}, the word length of $(A,k)\in C_2\wr C_n$ is equal to $|A|+\ell(A,k)$, where $\ell(A,k)$ is the length of the shortest path in the cycle $C_n$, going from $0$ to $k$ and containing $A$. From this it is clear that $diam(W^1_n)\leq 2n$.

\section{Comparison with similar inequalities}

Lower bounds of spectral nature on $c_{2}(X)$, can be traced back to \cite{LLR}. 
At least two other inequalities (see \cite{GN,NR}) linking the distortion, the $p$-spectral gap and other graph invariants have been published. In this section, we compare them to Theorem \ref{generalization}. We start with the Grigorchuk-Nowak inequality \cite{GN}.

\begin{Def}
Let $X$ be a finite metric space. Given $0<\epsilon <1$ define the constant $\rho_{\epsilon}(X)\in\left[0,1 \right]$, called the volume distribution, by the relation
$$
\rho_{\epsilon}(X)=\min \left\{\frac{diam(A)}{diam(X)} \ : \ A\subset X \; \mbox{such that $|A| \geq \epsilon |X|$}\right\} .
$$
\end{Def}

\begin{Thm}\label{GN}(\cite{GN} Theorem 3)
Let $X$ be a connected graph of degree bounded by $k$ and let $1\leq p< +\infty$. Then, for every $0<\epsilon <1$,
$$
\frac{(1-\epsilon)^{\frac{1}{p}}\rho_{\epsilon}(X)}{2^{\frac{1}{p}}} \ diam(X) \left( \frac{\lambda_{1}^{(p)}(X)}{k\ 2^{p-1}}\right)^{\frac{1}{p}} \leq c_{p}(X).
$$
\end{Thm}

It is easy to see that, when the graph satisfies $D(X)=diam(X)$ (this is the case for vertex-transitive graphs, by lemma \ref{Cayley}), then this result is weaker than our Theorem \ref{generalization}, since the factor $\frac{(1-\epsilon)^{\frac{1}{p}}\rho_{\epsilon}(X)}{2^{\frac{1}{p}}}$ is strictly smaller than $1$.

The second result, due to Newman-Rabinovich \cite{NR}, holds for $p=2$:

\begin{Prop}\label{NR}(\cite{NR} Proposition 3.2)
Let $X=(V,E)$ be a $k$-regular graph. Then,
$$
\sqrt{\frac{(|V|-1)\lambda_{1}^{(2)}(X)}{|V| \ k} \ avg(d^2)} \leq c_{2}(X),
$$
where $avg(d^2):=\frac{1}{|V| (|V|-1)}\sum_{x,y\in V} d(x,y)^2$.
\end{Prop}

In the following, we will compute the term $avg(d^2)$ for the cycle $C_{n}$ and for the hypercube $H_{d}$ in order to give explicitly the LHS term of the inequality due to Newman and Rabinovich. First, it is true that for a vertex-transitive graph $X=(V,E)$, we have
$$
\sum_{y,x\in V} d(x,y)^2 = |V| \sum_{j=1}^{diam(X)} j^2 |S(x_{0},j)|,
$$ 
where $x_{0}$ is an arbitrary point in $X$ and $S(x_{0},j)$ is the sphere of radius $j$, centered in $x_{0}$. By taking $n \geq 4$ and even, we clearly have
$$
\sum_{x,y\in C_{n}} d(x,y)^{2}=n\left(2\sum_{j=1}^{\frac{n}{2}-1} j^2 +\frac{n^2}{4}\right)= \frac{n^2(n^2+2)}{12}.
$$
Therefore, we get $\sqrt{\frac{n^2+2}{6}} \ \sin \frac{\pi}{n}$ as lower bound for $c_{2}(C_{n})$, which is strictly weaker than Corollary \ref{cycle}. On the other hand, for the hypercube $H_{d}$, by the same argument, we have
$$
avg(d^2)= \frac{1}{2^{d}(2^{d}-1)} \sum_{x,y\in H_{d}} d(x,y)^{2}= \frac{1}{2^{d}-1}\sum_{j=1}^{d}j^{2} {d \choose j}.
$$
Since $\sum_{j=1}^{d}j^{2} {d \choose j}<d^{2}2^{d-1}$ for $d\geq 2$, we conclude that Corollary \ref{hypercube} gives a better lower bound for $c_{2}(H_{d})$.

Finally, we mention for completeness a remarkable result, of a different nature, due to Linial, Magen and Naor \cite{LMN}: 

\begin{Thm}(\cite{LMN}, Theorem 1.3) There is a universal constant $C>0$ such that, for every $k$-regular graph $X$ with girth $g$:
$$c_2(X)\geq \frac{Cg}{\sqrt{\min\{g,\frac{k}{\lambda_1^{(2)}(X)}\}}}.$$
\end{Thm}

Observe however that, for the family $(H_d)_{d\geq 2}$ of hypercubes, the right-hand side of the inequality remains bounded, while $c_2(H_d)=\sqrt{d}$ by Corollary \ref{hypercube}.

\section{Bounds for compression}

In \cite{GK}, Guentner and Kaminker introduced a numerical quasi-isometry invariant called compression exponent to characterize how close to a quasi-isometry a coarse embedding can be. We recall the definitions. Let $(X,d)$ and $(Y,\delta)$ be two metric spaces. A map $F : X \rightarrow Y$ is a \textit{coarse embedding} of $X$ into $Y$ if there exist non-decreasing functions $\rho_{-},\rho_{+}: \mathbb{R} \rightarrow \mathbb{R}$ so that the following conditions hold:
\begin{enumerate}
	\item[(i)] $\rho_{-}(d(x,y)) \leq \delta(F(x),F(y)) \leq \rho_{+}(d(x,y))$, for all $x,y\in X$;
	\item[(ii)] $\lim_{t \rightarrow  +\infty}\rho_{-}(t)=+\infty$.
\end{enumerate}
$F$ is said to be \textit{large-scale Lipschitz} whenever $\rho_{+}$ can be taken of the form $t \mapsto Ct+D$, for some constants $C,D \geq 0$. In this case, we can define the \textit{compression} of a large-scale Lipschitz map $F$, denoted $R(F)$, as the supremum taken over all $\alpha\in \lbrack 0,1 \rbrack$ such that there exist $C,D \geq 0$ satisfying
$$
\delta(F(x),F(y)) \geq C \cdot  d(x,y)^{\alpha} -D, \ \ \forall x,y\in X.
$$
Then, the \textit{compression exponent} $\alpha_{Y}(X)$ corresponds to the supremum of $R(F)$ over all large-scale Lipschitz maps $F : X \rightarrow Y$. As target space, we will consider only $L^{p}=L^p([0,1])$. In this case, we write $\alpha_{p}(X)=\alpha_{L^{p}}(X)$. We will estimate the $L^{p}$ compression for a particular type of metric spaces obtained as disjoint union of finite graphs. \\
Let $(X_{n},d_{n})_{n}$ be a sequence of metric spaces. Let us denote by $d$ the metric on $X:=\sqcup_{n} X_{n}$, the disjoint union of the $X_{n}$, defined by
$$
d(x,y)=\left\{
\begin{array}{ll}
d_{n}(x,y), & \textrm{if}\ x,y\in X_{n};\\
\max \{diam(X_n), diam(X_{m})\}, & \textrm{if} \ x\in X_{n}, y\in X_{m}, \textrm{and} \ n\neq m .\\
\end{array}\right.
$$
The metric space $(X,d)$ is called the \emph{box space} of $(X_{n})_{n}$. The box space of a family of expanders was considered by Gromov to give the first example of metric spaces with bounded geometry that do not admit any coarse embedding into a Hilbert space (see \cite{Gro}).\\
In \cite{Aus}, Austin observed that it was possible to quantify the compression exponent of some infinite metric spaces containing arbitrarily large finite subsets with high distortions. Applying Austin's result to our setting,  we immediately get an estimate for the compressions of box spaces build from finite metric spaces.

\begin{Lem}\label{Austin}(\cite{Aus} Lemma 3.1.)
Let $V$ be a normed vector space and let $(X_{n},d_{n})_{n \geq 1}$ be a sequence of $1$-discrete, finite metric spaces satisfying the following 2 conditions: 
\begin{enumerate}
	\item[(i)] $\lim_{n \rightarrow + \infty}diam(X_{n})=+ \infty$;
	\item[(ii)] there exist $\eta >0$ and $K>0$ so that $c_{V}(X_{n})\geq K \cdot diam(X_{n})^{\eta}$, for all $n\geq 1$.
\end{enumerate}
Then, denoting $X$ the box space of $(X_{n})_{n\geq 1}$, we have $\alpha_{V}(X) \leq 1-\eta$.
\end{Lem}

\begin{Cor}
Using the notations of Proposition \ref{slnp}, $\alpha_{2}(\sqcup_{n\geq 2} Y_{n})\leq \frac{3}{4}$.
\end{Cor}

{\bf Proof:} Since $c_{2}(Y_{n})=\Omega(diam(Y_{n})^{\frac{1}{4}})$ (see the proof of Proposition \ref{slnp}), the result follows immediately from lemma \ref{Austin}. 
\hfill $\square$

\begin{Cor}
$\alpha_{2}(\sqcup_{d\geq 2} H_{d})=\frac{1}{2}$.
\end{Cor}

{\bf Proof:} The first inequality $\alpha_{2}(\sqcup_{d\geq 1} H_{d}) \leq\frac{1}{2}$ is obvious from Corollary \ref{hypercube} and lemma \ref{Austin}. To prove the converse inequality, we simply remark that $\sqcup_{d\geq 2} H_{d}$ admits a bi-Lipschitz embedding into $\ell^{1}(\mathbb{N})$ and the conclusion follows from Proposition 2.6 in \cite{GK}. Indeed, let $\phi_{k} : \sqcup_{d\geq 2} H_{d} \rightarrow \mathbb{R}^{k+1}: x \mapsto \phi_{k}(x)$, where
$$
\phi_{k}(x)=\left\{
\begin{array}{ll}
(diam(H_{k}),x), & \textrm{if}\ x\in H_{k};\\
0, & \textrm{otherwsie}.\\
\end{array}\right.
$$
Then, $\phi$, defined by $\phi=\oplus_{d\geq 2} \phi_{d}$, where the target space is endowed with the $\ell^{1}$ norm, is a bi-Lipschitz embedding.
\hfill $\square$

\medskip
The interest of this particular box space is due to the fact that it is the first explicit example of a metric space without Yu's property A which is coarsely embeddable into a Hilbert space (see \cite{Now}).

\bigskip

Authors addresses:
\medskip

\noindent
Institut de Math\'ematiques - Unimail\\
11 Rue Emile Argand\\
CH-2000 Neuch\^atel\\
Switzerland

\begin{verbatim}pierre-nicolas.jolissaint@unine.ch; alain.valette@unine.ch
\end{verbatim}

\end{document}